# On the relation of truncation and approximation errors for the set of solutions obtained by different numerical methods

Alekseev A.K., Bondarev A.E.
Keldysh Institute of Applied Mathematics RAS
Moscow, Russia
e-mail: aleksey.k.alekseev@gmail.com, bond@keldysh.ru

*Abstract*

The truncation and approximation errors for the set of numerical solutions computed by methods based on the algorithms of different structure are calculated and analyzed for the case of the two-dimensional steady inviscid compressible flow. The truncation errors are calculated using a special postprocessor, while the approximation errors are obtained by the comparison of the numerical solution and the analytic one. The extent of the independence of errors for the numerical solutions may be estimated via the Pearson correlation coefficient that may be geometrically expressed by the angle between errors. Due to this reason, the angles between the approximation errors are computed and related with the corresponding angles between the truncation errors. The angles between the approximation errors are found to be far from zero that enables a posteriori estimation of the error norm. The analysis of the distances between these solutions provides another approach to the estimation of the error. The comparison of the error norms, obtained by these two procedures, is provided that demonstrates the acceptable values of their effectivity indices. The results of the approximation error norm estimation for the supersonic flows, containing shock waves, are presented. The measure concentration phenomenon and the algorithmic randomness give some insights into these results.

**1. Introduction**

The estimation of an approximation error is of the high current interest at the numerical solutions of Partial Differentiation Equations [1-3]. The significant number of the computationally efficient error estimators is available in the domain of elliptic equations [1] that is caused by the high regularity of corresponding solutions. At present, the progress of the error estimation for Computational Fluid Dynamics (CFD) problems is more limited due to the presence of discontinuities that may spontaneously arise in the flow domain. This circumstance enforces additional search for the construction of new error estimators and an investigation of the properties of truncation and approximation errors. Several new approaches apply the set of independent numerical solutions [4-8] for the elaboration of error norm estimators. Herein, we consider the relation of these methods and the measure concentration effect [9-16] that may shed an addition light on the numerical results.

The measure concentration effect relates the probability and the geometry of the domains in multidimensional spaces (hyperspheres, hypercubes etc.). For example, the Gaussian distribution is defined by the concentration of the mass in a hypersphere near its surface. The measure concentration provides certain unexpected but powerful feasibilities for analysis that can be called the "blessing of dimensionality" [13].

---

The measure concentration effect is of the great importance for the statistic thermodynamics as the system of many variables. From this standpoint, the numerical solutions (grid functions), appearing at a discretization of multidimensional partial differential equations, present an interesting analogy, since the dimension of the grid function space $R^N$ (number of grid nodes multiplied by the number of the flow variables) is great. It may take the values $N \sim 10^5 \div 10^6$ and higher for standard CFD problems. Above the high dimensionality and the shape of the grid functions domain (for example, hypersphere), the arbitrary (random) choice of grid functions, is crucial for the applicability of the measure concentration effect.

We consider a system of CFD equations in the operator form $A\widetilde{u} = f$ and corresponding discrete operator $A_h u_h = f_h$ that approximates the system on some grid. In further presentation we consider the numerical solution $u_h$ to be a grid function (vector $u_h \in R^N$), $\widetilde{u}_h \in R^N$ to be the projection of true solution onto the grid, $\Delta \widetilde{u}_h = u_h - \widetilde{u}_h \in R^N$ to be the true approximation error, $\Delta u_h \in R^N$ to be an estimate of the approximation error. The truncation error $\delta u \in R^N$ may be assesed using the Taylor series expansion of the grid function and substitution of result to the discrete algorithm. We consider the set of $K$ independent numerical methods $A_h^k u^{(k)} = f_h (k=1...K)$, the numerical solutions $u^{(k)}$, obtained by these methods, and corresponding truncation $\delta u^{(k)}$ and approximation $\Delta u^{(k)}$ errors. By the "independent" algorithms we mean structurally different algorithms, for example, methods of the different nominal order of approximation. We understand the "independent" numerical solutions in the corresponding sense.

The properties of the truncation and approximation errors on the ensemble of independent solutions obtained on the same grid are the main topic of the current paper.

The numerical solutions of CFD problems obtained by independent algorithms approximate the exact solution and, in the ideal case, should be close in some norm. Thus, the set of these grid functions should be highly correlated and be out our attention.

The truncation and the approximation errors are also correlated with the flowfield, yet they significantly depend on the algorithm of computations that may provide some decorrelation. For example, the truncation error for calculations, performed by algorithms of the different order of the approximation, is composed by the derivatives of different leading (least) orders. From this viewpoint, these truncation errors are linearly independent for nontrivial solutions. The analysis of the truncation errors' structure provides some arguments that support the algorithmic randomness of the corresponding grid functions and engender prospects for their decorrelation.

The approximation error is caused by the truncation error transformation performed by the solution operator. The artificial viscosity or flux/slope limiters engender an additional error connected with the discontinuities in the flowfield and, by this reason, highly correlated even for algorithms of the different structures. So, the correlation of the approximation error for independent solutions may be significant and needs for an investigation.

The aim of the current paper is to study the correlation of truncation and the approximation errors from the prospects for a posteriori error (approximation error norm) estimation using the ensemble of the numerical solutions obtained by the independent methods. In contrast to the deterministic approaches [1-5] and to stochastical one [17] we consider the line of attack that is based on the algorithmic randomness [18-20].

The paper is organized as follows. In Section 2 we consider some properties of the truncation and approximation errors. The algorithmic randomness of the truncation error is discussed in line with the orthogonality of random vectors in high dimensions. In Section 3 we consider the geometry of the high-dimensional spaces from the standpoint of the approximation error estimation with account of deviations from the pure orthogonality. The results of numerical tests are presented in Section 4 for inviscid supersonic flows computed by the set of independent methods. The correlation of the approximation and truncation errors is analyzed. The efficiency indices of the approximation error norm estimation are presented and discussed in Section 5.

## 2. Stochastic features of the truncation and approximation errors

Let we have the set of $K$ numerical algorithms $A_h^k u^{(k)} = f_h$ that approximate the system $A\widetilde{u} = f$ on some fixed grid. The truncation error (source term of differential approximation [21]) may be formally expressed as

$$\delta u^{(k)} = \sum_{m=j_k}^{\infty} C_m^{(k)} h^m \frac{\partial^{m+1} \widetilde{u}}{\partial x^{m+1}} \quad (1)$$

($j_k$ is the approximation order of $k-th$ algorithm). Unfortunately, this expression contains the infinite number of members that restricts its computability by terms of some (usually, minor) order.

The finite (Lagrange) form of truncation error has the appearance

$$\delta u^{(k)}(\alpha_n) = C^{(k)} h^{j_k} \cdot \partial^{j_k+1} \widetilde{u}(x_n + \alpha_n h) / \partial x^{j_k+1} \quad (2)$$

($n$ is the grid node number) [22,23]. It contains an unknown parameter $\alpha_n \in [0,1]$. Numerical tests performed for the heat conduction equation [22] demonstrate the parameter $\alpha_n$ to obey the probability distribution density $P(\alpha_n)$ of rather universal form (although, non Gaussian). So, the truncation error reveals certain probabilistic features despite its deterministic origin.

On the other hand, several papers ([17] for example) assume the discretization error $\Delta u^{(k)} = (u^{(k)} - \widetilde{u}_h)$ to be a random normally distributed value that, at first glance, appears to be unexpected.

So, the relation of the random and deterministic properties of the truncation and approximation errors is not trivial. This circumstance may be clarified using the measure concentration theory relating the probability and the geometry of multidimensional spaces [10-16]. The dimensionality of the space, the geometry of the domains (compact subsets) containing errors (hypersphere, hypercube, etc) and the random choice of vectors are the main factors that determine the considered properties of errors. Herein we discuss the application of the measure concentration theory and the algorithmic randomness to the estimation of the discretization error norm.

### 2.1. The algorithmic randomness

The calculations performed by the structurally different numerical algorithms may provide the low correlation of truncation errors $\delta u^{(k)}$ for independent solutions, since $\delta u^{(k)}$ depends not only on the solution $u^{(k)}$, but on the choice of the algorithm (the order of the minor derivatives and corresponding coefficients for schemes of different order of approximation, for example). So, herein, we consider the conjectured independence (decorrelation) of truncation errors from the standpoint of the algorithmic randomness. By definition ([18,19]) the vector $\delta u^{(k)}$ is algorithmic random, if it is shorter any code that may reproduce it. The truncation error series (Eq. 1) starts from the derivatives of the different orders for algorithms of different order of approximation. By this reason, such series are linearly independent (if exclude some trivial solutions). The transformation $u^{(k)} \to \delta u^{(k)}$ contains the infinite sum of derivatives. From this standpoint, it requires an algorithm of the infinite length for the calculation of the finite vector $\delta u^{(k)}$. This circumstance may be interpreted as the evidence of the algorithmic randomness ("the string is shorter the program"). If the algorithmic randomness is considered as the random choice with restrictions of zero measure [20] (rejection of too "regular" vectors), then the ensemble of calculations performed by the algorithms with different approximation orders corresponds to the set of randomly chosen $\delta u^{(k)}$. So, the conjecture of $\delta u^{(k)}$ algorithmic (Kolmogorov) randomness seems to be justified and should be checked in numerical tests.

*2.2 The orthogonality in spaces of the great dimension*

It is known [10-14] that two vectors $v^{(1)} \in R^N$ and $v^{(2)} \in R^N$ randomly selected from two hyperspheres in spaces of the great dimension $N$ are orthogonal with the "*probability one*" ($P$ is the probability, $(\cdot,\cdot)$ is the scalar product):

$$P\{(v^{(1)}/\|v^{(1)}\|, v^{(2)}/\|v^{(2)}\|) > \delta\} < \sqrt{\pi/2}\, e^{-\delta^2 N/2}. \tag{3}$$

Herein and below we consider the scalar product and the norm in the space $L^2$. In general, the high dimensional space has a huge number of orthogonal vectors with the exponential dependence on the space dimension. For another example, all diagonal vectors of the hypercube are orthogonal to all the axes. All these objects are orthogonal with the "*probability one*".

Thus, the deterministic objects in high dimensions manifest some features of the probabilistic nature, such as the normal distribution or some properties that appear with the "*probability one*". This may explain results of papers [17, 22] related with the unexpected statistical properties of the truncation and approximation errors. The more information on the relation of the geometry and probability may be found in [15, 16].

Both the truncation error $\delta u_h^{(k)} \in R^N$ and the approximation error $\Delta u^{(k)} = u^{(k)} - \widetilde{u}_h \in R^N$ belong to the spaces of the grid functions of very high dimension. Due to the measure concentration effect, this may cause the orthogonality (in sense of "probability one"). The truncation $\delta u^{(k)}$ and approximation errors $\Delta u^{(k)}$ origin at zero and belong to some unknown parts of a hypersphere that may cause deviations from the exact orthogonality. The geometry of domains that contain $\delta u^{(k)}$ and $\Delta u^{(k)}$ is unknown and, so, is of the current interest from prospects of the applicability to a posteriori error estimation. From this standpoint we study the values of scalar products $(\delta u^{(k)}, \delta u^{(i)})$ and $(\Delta u^{(k)}, \Delta u^{(i)})$ and the corresponding angles between the vectors in this paper.

**3. The estimation of the approximation error norm**

One may easy see the difference of two numerical solutions to be equal to the difference between their errors $u^{(1)} - u^{(2)} = u^{(1)} - \widetilde{u}_h - u^{(2)} + \widetilde{u}_h = \Delta u^{(1)} - \Delta u^{(2)}$. This relation contains certain information on the unknown errors that we want to use.

Let's assume these errors $\Delta u^{(1)}, \Delta u^{(2)}$ to belong hyperspheres with the center at zero point. If two errors are randomly chosen and belong to the high-dimension space, the distance $d_{1,2} = \|u^{(1)} - u^{(2)}\| = \|\Delta u^{(1)} - \Delta u^{(2)}\|$ between numerical solutions $u^{(1)} \in R^N$ и $u^{(2)} \in R^N$ "*with probability one*" is greater the distance from the exact solution to the numerical one (the approximation error), since $\Delta u^{(k)}$ "*with probability one*" are orthogonal and the hypotenuse (distance between solutions) is greater any leg

$$d_{1,2} = \|\Delta u^{(1)} - \Delta u^{(2)}\| \geq \|\Delta u^{(k)}\| = \|u^{(k)} - \widetilde{u}_h\|, k = 1,2. \tag{4}$$

So, the "blessing of dimensionality" [13] may provide a posteriori error estimation on the ensemble of independent solutions.

The expression (4) resembles the results obtained by the triangle inequality [4]

$$\left\| u^{(2)} - \tilde{u}_h \right\| \leq d_{1,2}. \tag{5}$$

The expression (5) is based on a priori known relation of error norms for two numerical solutions, does not assume the error orthogonality, and provides an upper estimate only for the more precise solution $u^{(2)}$.

If the errors are strictly orthogonal, the famous hypercircle method [24] provides the precise estimate

$$\left\| (u^{(1)} + u^{(2)})/2 - \tilde{u}_h \right\| = d_{1,2}/2. \tag{6}$$

Unfortunately, the exact orthogonality is not observed in the numerical tests and there exists the need for more robust estimates that are considered herein.

The error estimate by Eq. (4) for $k-th$ solution may be naturally extended on the ensemble of $K$ solutions as follows:

$$d_{k,\max} = \max_i \left\| u^{(k)} - u^{(i)} \right\|, (i=1,K), \tag{7}$$

$$d_{k,\max} \geq \left\| \Delta u^{(k)} \right\|. \tag{8}$$

Using the ensemble width (the maximum distance between solutions on the ensemble)

$$d_{\max} = \max_{i,k} \left\| u^{(k)} - u^{(i)} \right\|, \tag{9}$$

we may obtain more pessimistic estimate, suitable for all considered solutions:

$$d_{\max} \geq \left\| \Delta u^{(k)} \right\|. \tag{10}$$

These estimates (Eqs. (8) and (10)) seem to become more reliable as the number of the ensemble elements expands.

*3.1 The bounds of the error estimation for nonorthogonal errors*

The estimate by Eq. (4) is based on the assumption of approximation error $\Delta u^{(k)}$ orthogonality. Unfortunately, the approximation error $\Delta u^{(k)}$ is highly correlated near discontinuities, since it is a wave with the positive and negative wings for a monotonic solution. The artificial viscosity enforces the monotonicity in vicinity of discontinuities for modern numerical methods. Thus, $\Delta u^{(k)}$ are not exactly orthogonal and the estimate (4) may be violated.

Fortunately, the inequality (4) is enough robust and valid at the significant deviation from the orthogonality. The value of the angle between errors $\alpha = \arccos\theta$ ($\theta = \frac{(\Delta u^{(1)}, \Delta u^{(2)})}{\left\|\Delta u^{(1)}\right\| \cdot \left\|\Delta u^{(2)}\right\|}$ is the Pearson correlation coefficient) is of the interest from the viewpoint of the relation between the independent and correlated components of the error. One may see from the elementary geometric considerations the inequality (4) to be valid for enough great angles $\alpha \geq \alpha_*$. The minimum value $\alpha_* = 60^o$ occurs if the norms of errors are equal.

*3.2 The estimation of the error's norm with account of the angle between vectors of errors*

In this Section we show the computable value $\|\Delta u^{(1)} - \Delta u^{(2)}\| = \|u^{(1)} - u^{(2)}\|$ to provide the estimate for $\|\Delta u^{(k)}\|$ if the angle $\alpha$ between errors is known. Let's analyze the Fig. 1, where the norms of unknown errors are marked as $\|\Delta u^{(1)}\| = r_1$ and $\|\Delta u^{(2)}\| = r_2$. We assume $r_1 < r_2$ without a loss of generality. The distance $d_{1,2} = \|\Delta u^{(1)} - \Delta u^{(2)}\| = \|u^{(1)} - u^{(2)}\|$ and the angle $\alpha$ between errors are known. We form the auxiliary isosceles triangle with sides $r_1$. The estimate $d_0 < d_{1,2} \cdot \cos(\varphi)$ may be done for the third side of this triangle. Then the following relation is valid

$$r_1 = \frac{d_0}{2\sin(\alpha/2)} < \frac{d_{1,2} \cdot \cos(\varphi)}{2\sin(\alpha/2)}. \tag{11}$$

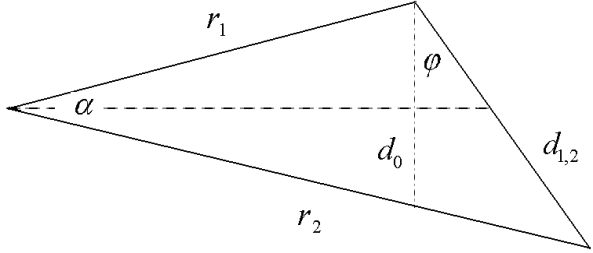

Fig. 1. The relation of the errors, the distance between solutions and the angle between errors.

The second side equals $r_2 = r_1 + d_{1,2} \dfrac{\sin(\varphi)}{\cos(\pi - \pi/2 - \pi/2 + \alpha/2)} = r_1 + d_{1,2} \dfrac{\sin(\varphi)}{\sin(\alpha/2)}$.

By account of Eq. (11) one obtains

$$r_2 < \frac{d_{1,2} \cdot \cos(\varphi)}{2\sin(\alpha/2)} + \frac{d_{1,2} \sin(\varphi)}{\sin(\alpha/2)} = \frac{d_{1,2}}{2\sin(\alpha/2)}(\cos(\varphi) + 2\sin(\varphi)). \tag{12}$$

The maximum for the right hand side of this expression occurs at $\sin(\varphi) = 2\cos(\varphi)$, so $tg(\varphi) = 2$, and $\varphi_{max} = \arctan(2)$ ($\varphi_{max} \approx 63.43$). In result, the following relation is valid:

$$\|\Delta u^{(k)}\| < \frac{\|u^{(1)} - u^{(2)}\|}{2\sin(\alpha/2)}(\cos(\varphi_{max}) + 2\sin(\varphi_{max})), (k = 1,2), \varphi_{max} = \arctan(2). \tag{13}$$

This expression may be approximated by

$$\|\Delta u^{(k)}\| < 1.1 \cdot \frac{\|u^{(1)} - u^{(2)}\|}{\sin(\alpha/2)}, (k = 1,2) \tag{14}$$

that enables the estimation of error norm using the value of the angle $\alpha$ between approximation errors and the distance between the numerical solutions $d_{1,2} = \|u^{(1)} - u^{(2)}\|$.

The information regarding the angle between approximation errors $\alpha$ may be obtained from the angle $\beta$ between truncation errors and the empirical relation $\alpha(\beta)$ that is considered below.

## 4. Numerical results

The above analysis is verified by the set of numerical tests. The numerical solutions for the following system of equations are used in the tests:

$$\frac{\partial \rho}{\partial t} + \frac{\partial (\rho U^j)}{\partial x^j} = 0, \tag{15}$$

$$\frac{\partial (\rho U^i)}{\partial t} + \frac{\partial (\rho U^j U^i + P\delta_{ij})}{\partial x^j} = 0, \tag{16}$$

$$\frac{\partial (\rho E)}{\partial t} + \frac{\partial (\rho U^j h_0)}{\partial x^j} = 0. \tag{17}$$

This system of equations describes flows of inviscid compressible fluid. The system is two-dimensional, $i,j = 1,2$, $U^1 = U, U^2 = V$ are the velocity components, $h_0 = (U^2 + V^2)/2 + h$, $h = \frac{\gamma}{\gamma-1}\frac{P}{\rho} = \gamma e$, $e = \frac{RT}{\gamma-1}$, $E = \left(e + \frac{1}{2}(U^2 + V^2)\right)$ are enthalpy and energies, $P = \rho RT$ is the state equation, $\gamma = C_p/C_v$ is the specific heats relation.

The flow patterns (Fig. 2 and Fig. 3) occurring at the shock waves interaction (Edney-I and Edney-VI flow structures, [25]) are considered, since they have the analytical solutions. The flows are steady and computed by the relaxation over the time from the uniform initial conditions. The projections of analytical solutions to the grid $\tilde{u}_h$ are considered as the "exact" solutions. The "exact" discretization error $\Delta u^{(k)} = u^{(k)} - \tilde{u}_h$ was calculated by subtraction of the numerical solution from the projection of the analytical one onto the computational grid.

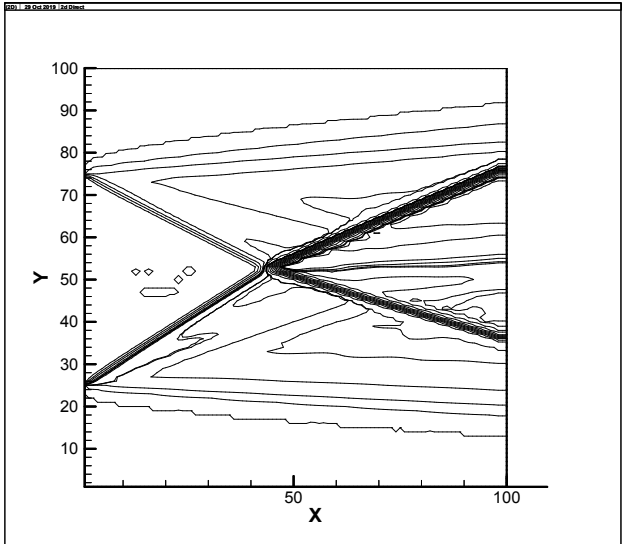 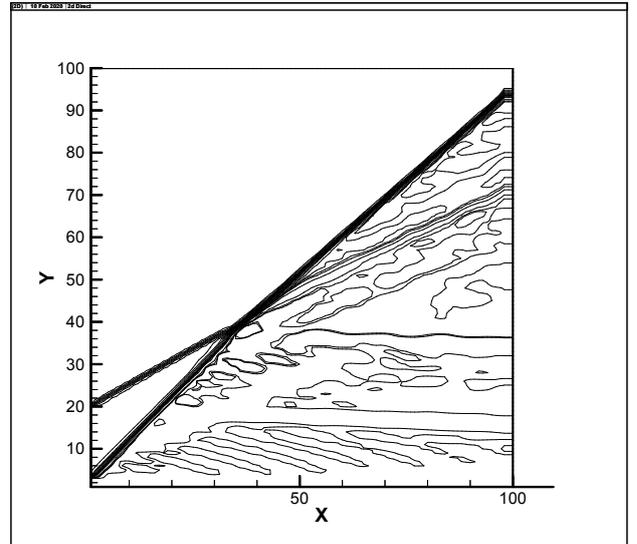

Fig. 2. Density isolines (Edney I)  Fig. 3. Density isolines (Edney VI)

Fig. 2 provides the density distribution for Edney-I flow pattern ($M = 4$, flow deflection angles $\alpha_1 = 20^o$ (lower) and $\alpha_2 = 15^o$ (upper)). Fig. 3 provides the density distribution for Edney-VI ($M = 3.5$, consequent angles $\alpha_1 = 15^o$ and $\alpha_2 = 25^o$). Both flowfields are determined by the crossing shock waves and contact discontinuities.

We analyze the set of 13 numerical solutions obtained on the same grid by the numerical methods [26-40] covering the range of approximation order from one to five. The tests for Edney I and Edney VI flows are performed using uniform grids of $100 \times 100$ and $400 \times 400$ nodes.

*4.1 The estimation of truncation error*

The computation of the approximation error is the elementary operation if the analytic solution is available. The complete computation of the exact truncation error is not feasible due to the presence of the infinity series or unknown Lagrange coefficients. However, the estimation of the leading part of the truncation error may be successfully performed. There exist sophisticated methods for the truncation error estimation [41,42] related with the consequent interpolation and differentiation of a grid function. Herein, the truncation error $\delta u^{(k)}$ was estimated using postprocessor described by [43]. It is the simplest approach that is sketched below. Let's have the flowfield computed by some (in general, unknown) finite difference (volume) algorithm $u_h^{(k)} = (A_h^{(k)})^{-1} f_h$. We seek for the numerical estimate of the truncation error $\delta u^{(k)}$. Let the finite differences be expanded in the Taylor series and the differential approximation $A u_h^{(k)} = f_h + \delta u_h^{(k)}$ be obtained. At the next step, we act by certain stencil of the higher approximation order $A_h^{high}$ (postprocessor) on the computed flowfield (grid function) and obtain

$$A_h^{high} u_h^{(k)} = A_h^{high} (\widetilde{u}_h + \Delta u_h^{(k)}) =$$
$$= A_h^{high} (u_h^{(high)} - \Delta u_h^{(high)} + \Delta u_h^{(k)}) = f_h - \delta u_h^{(high)} + A_h^{high} (A_h^{(k)})^{-1} \delta u_h^{(k)}.$$

The discrepancy $A_h^{high} u_h^{(k)} - f_h$ estimated by the action of the high order stencil on the low order solution (grid function) has the appearance

$$\eta_h = A_h^{high} u_h^{(k)} - f_h = -\delta u_h^{(high)} + A_h^{high} (A_h^{(k)})^{-1} \delta u_h^{(k)} \approx -\delta u_h^{(high)} + \delta u_h^{(k)}. \qquad (18)$$

Herein we assume $A_h^{high} (A_h^{(k)})^{-1} \approx E$ ($E$ is the unite matrix) since both used operators approximate the same $A$. In the result, the misfit $\eta_h$ contains truncation errors both of the investigated algorithm and of the postprocessor. We consider this value as the estimate of the truncation error by neglecting the postprocessor error, since it has the greater order over the grid step. The high order stencil operates herein only locally (for the precomputed flowfield analysis), so the stability and monotonicity properties may be not considered.

The postprocessor, based on the stencil of sixth approximation order is used:

$$\frac{-f_{k+3}^n + 9 f_{k+2}^n - 45 f_{k+1}^n + 45 f_{k-1}^n - 9 f_{k-2}^n + f_{k-3}^n}{60 h_k}. \qquad (19)$$

*4.2 The correlation of truncation and approximation errors on the set of solutions*

The angles between truncation errors $\delta u^{(k)}$ are calculated for the ensemble of $K=13$ solutions obtained using 10 numerical algorithms of different structure and properties [26-40] and several variants of the artificial viscosity. The angles between $\Delta u^{(k)}$ are also calculated for this ensemble. The numerical solutions are obtained using following methods:

First order algorithm by Courant-Isaacson-Rees) [26] implemented in accordance with [27],
Second order MUSCL [28] based algorithm that uses approximate Riemann solver by [29],
Second order MUSCL [28] based algorithm that uses approximate Riemann solver by [40],
Second order algorithm based on the relaxation approach by [33],
Second order algorithm based on the MacCormack [34] scheme: without the artificial viscosity, with the second order artificial viscosity with coefficients $\mu=0.01$ and $\mu=0.002$, with the fourth order artificial viscosity ($\mu=0.01$),
Second order algorithm based on the "two step Lax- Wendroff" method [35,36] with the artificial viscosity of second order ($\mu=0.01$),
Third order algorithm based on the modification of Chakravarthy-Osher method [30,31],
Third order algorithm WENO [39],
Fourth order algorithm by [32],
Fifth order algorithm WENO [37, 38].

The angles between truncation errors $\beta_{km}=\arccos((\delta u^{(k)}\cdot\delta u^{(m)})/(\|\delta u^{(k)}\|\cdot\|\delta u^{(m)}\|))$ and the angles between approximation errors $\alpha_{km}=\arccos((\Delta u^{(k)}\cdot\Delta u^{(m)})/(\|\Delta u^{(k)}\|\cdot\|\Delta u^{(m)}\|))$ are calculated. The ensemble of $K=13$ solutions, which is used in tests, provides $K\cdot(K-1)/2=78$ values of angles $\beta_{km}$ and $\alpha_{km}$. These angles are computed for two flow patterns and two grids, the results are provided in Fig. 4 as $\alpha_{km}(\beta_{km})$. The averaged (over the ensemble) angle between the truncation errors is in the range $\bar{\beta}=58\div64$ for the considered flow patterns and grids. The averaged angle between the approximation errors is in the range $\bar{\alpha}=30\div44$.

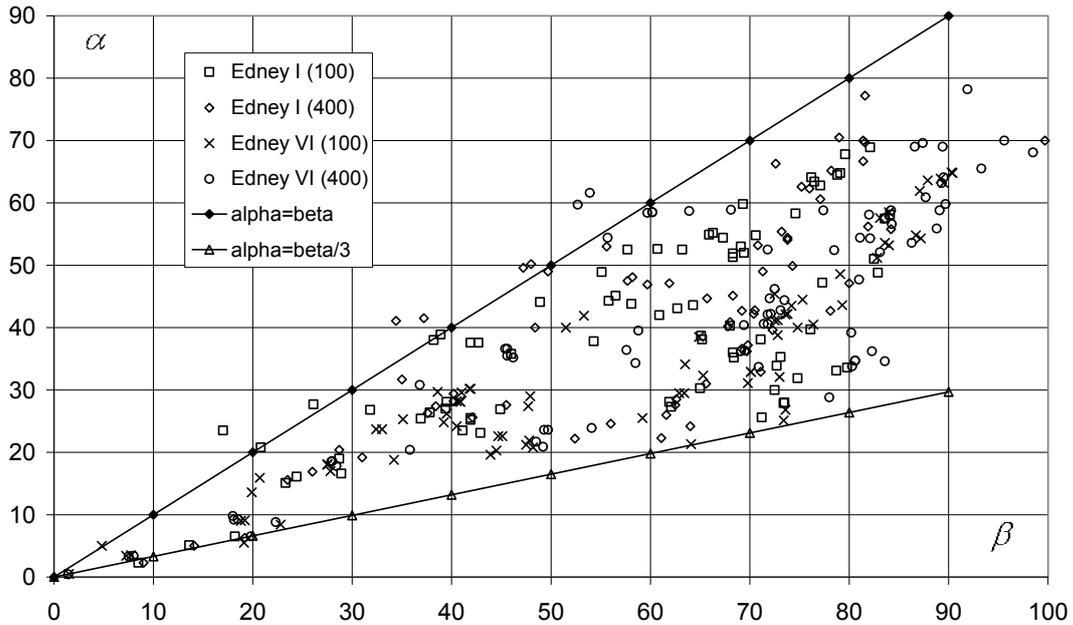

Fig. 4. The angle between approximation errors in dependence on the angle between the truncation errors.

The angle between approximation errors $\alpha_{km}$ is related with the error between truncation errors $\beta_{km}$ by the expression:

$$\beta_{12} = \arccos\left[\frac{(\delta u^{(1)}, \delta u^{(2)})}{(\delta u^{(1)}, \delta u^{(1)})^{1/2} \cdot (\delta u^{(2)}, \delta u^{(2)})^{1/2}}\right] = \arccos\left(\frac{(A_h^{(1)}\Delta u^{(1)}, A_h^{(2)}\Delta u^{(2)})}{(A_h^{(1)}\Delta u^{(1)}, A_h^{(1)}\Delta u^{(1)})^{1/2} \cdot (A_h^{(2)}\Delta u^{(2)}, A_h^{(2)}\Delta u^{(2)})^{1/2}}\right) =$$

$$= \arccos\left(\frac{(\Delta u^{(1)}, A_h^{(1)*} A_h^{(2)} \Delta u^{(2)})}{(\Delta u^{(1)}, A_h^{(1)*} A_h^{(1)} \Delta u^{(1)})^{1/2} \cdot (\Delta u^{(2)}, A_h^{(2)*} A_h^{(2)} \Delta u^{(2)})^{1/2}}\right).$$

The properties of the operators $(A_h^{(k)})^* A_h^{(m)}$ determine the relation between $\alpha_{km}$ and $\beta_{km}$ (if $(A_h^{(k)})^* A_h^{(m)} = E$ then $\alpha_{km} = \beta_{km}$). By this reason (and with account of $A_h^{(m)}$ to be the approximation of the exact operator $A$) we assume the existence of certain unknown dependence $\alpha_{km}(\beta_{km})$ that is determined by the properties of the primal operator (more precisely $A^*A$). Fig. 4 demonstrates some results of the search for this dependence. The inequality $\beta_{km} > \alpha_{km}$ holds almost in all tests. The relation of $\alpha_{km}$ and $\beta_{km}$ may be expressed as $\alpha_{km} \geq \beta_{km}/3$. This inequality may be used for the error estimation.

*4.3 The comparative quality of considered error estimators*

The distances between solutions $d_{km} = \|u^{(k)} - u^{(m)}\| = \|\Delta u^{(k)} - \Delta u^{(m)}\|$ may be used for the error norm $\|\Delta u^{(k)}\|_{L_2}$ estimations in accordance with Eqs. (4), (8), (10) or (14).

In order to describe the quality of a posteriori error estimation we apply the effectivity index [1], which is equal to the relation of the norm of the error estimate to the norm of the true error

$$I_{eff} = \frac{\|\Delta u\|}{\|\Delta \widetilde{u}\|}. \tag{20}$$

A reasonable error estimate should be greater the error on the one hand ($I_{eff} \geq 1$) and be not too pessimistic on the other. The relation $I_{eff} < 3$ is stated by [1] that is specified for the finite-element applications and smooth solutions. For the CFD applications with shock waves one may expect this range to be wider. Corresponding values should be established in numerical tests.

The above discussed error norm estimators demonstrate different values of effectivity indices for the considered numerical tests.

A posteriori error estimation by the distance between two solutions (Equation (4) ( $d_{1,2} = \|u^{(1)} - u^{(2)}\| \geq \|\Delta u^{(k)}\|, k = 1,2$)) demonstrates $I_{eff} \sim 0.04 \div 1.5$ and, so, may be not reliable.

By increasing the number of the ensemble elements one may obtain the error norm estimates of more acceptable quality for a great enough size of the ensemble.

If one uses the expression (8) for the error norm estimation, the effectivity index acquires the form

$$I_{eff,k} = \frac{\max_m (d_{km})}{\|\Delta u^{(k)}\|} = \frac{\max_m (d_{km})}{\|u^{(k)} - \widetilde{u}\|}. \tag{21}$$

In above described numerical tests, the effectivity index based on (21) is in the range $I_{eff,k} \sim 1.1 \div 1.5$.

The effectivity index values based on the ensemble width (10)

$$I_{eff} = \max_{i,j}(d_{ij})/\|\Delta u^{(i)}\| \tag{22}$$

are in the range $I_{eff} \sim 1.1 \div 1.8$.

The value of effectivity index for (8) and (10) asymptotically approaches to results obtained by the triangle inequality [4, 5] as the number of solutions increases.

The estimation of the angle between truncation errors and the employment of inequality $\alpha_{km} \geq \beta_{km}/3$ enable the application of the Expression (14) $(1.1 \cdot \frac{\|u^{(1)} - u^{(2)}\|}{\sin(\alpha/2)} > \|\Delta u^{(k)}\|, (k=1,2))$.

The value of the effectivity index of this estimate is in the range $I_{eff} \sim 0.9 \div 4.5$. This approach provides the results of acceptable quality with the minimal computational burden (using only two solutions).

## 5. Discussion

The application of the ensemble of independent solutions for a posteriori error estimation requires moderate computational resources that are equivalent to several runs of the code on the same grid. It is much more inexpensive if compare with the Richardson extrapolation that requires several consequent refinements of grids. The need for the set of different codes based on the independent algorithms is more restrictive. Fortunately, at present, there exists a significant number of open access codes. For example, five different algorithms from OpenFoam are used in [6,7].

## 6. Conclusion

The numerical tests demonstrated that the ensemble of numerical solutions for Euler equations obtained by independent algorithms enables the estimation of the approximation error norm using the distances between solutions. The increase of the number of samples in the ensemble improves the reliability of the error estimation.

The analysis of tests demonstrated the angles between the truncation errors $\beta_{km}$, averaged over the ensemble, to be in the range $58^o \div 64^o$. The averaged angles between the approximation errors $\alpha_{km}$ are in the range $30^o \div 44^o$. The scatter of results corresponds different flow modes and grid sizes. The observed correlation of the approximation and truncation angles ($\alpha_{km} \geq \beta_{km}/3$) enables the error norm estimation using only two independent solutions.

The systematically great angles between the truncation errors $\beta_{km}$ may be attributed to the measure concentration phenomenon and the algorithmic randomness of the truncation errors for independent numerical algorithms.